 \documentclass[12pt,a4paper]{amsart}
  \usepackage{latexsym} 
  \usepackage[all]{xy}
  \usepackage{amsfonts} 
  \usepackage{amsthm} 
  \usepackage{amsmath} 
  \usepackage{amssymb}
  \usepackage{pifont}  
  
  \usepackage[T1]{fontenc}
\normalfont

  \def\<{{\langle}} 
  \def\>{{\rangle}}

  \def\eps{\varepsilon}

  \def\note#1{{}}

 \def\graph#1{{\rm \sf{graph}}(#1)} 

  \def\note#1{}

  \def\C{\mathbb{C}}

  \def\beq{\begin{equation}} 
  \def\eeq{\end{equation}}

  \def\id{\mathrm{id}} 
  \def\im{{\rm Im}}

  \def\ot{{\otimes}}

   \def\roR{\varrho^{R}} 
  \def\Rro{{}^{R}\!\varrho} 
  \def\Sro{{}^{S}\!\varrho}

     \def\1{\mathbf{1}}
     \def\coten#1{\Box_{#1}}

\def\k{\Bbbk}

  \newcounter{zlist}

  \newcounter{blist} 
  \newenvironment{blist}{\begin{list}{(\alph{blist})}{ 
  \usecounter{blist}\leftmargin2.5em\labelwidth2em\labelsep0.5em 
  \topsep0.6ex 
  \parsep0.3ex plus0.2ex minus0.1ex}}{\end{list}} 

  \newcounter{rlist}

   \newcounter{alist}

\def\stac#1{\raise-.2cm\hbox{$\stackrel{\displaystyle\otimes}{\scriptscriptstyle{#1}}$}}

\def\cten#1{\raise-.2cm\hbox{$\stackrel{\displaystyle\widehat{\otimes}}
{\scriptscriptstyle{#1}}$}}

  \def\Label#1{\label{#1}\ifmmode\llap{[#1] }\else 
  \marginpar{\smash{\hbox{\tiny [#1]}}}\fi} 
  \def\Label{\label}

  \theoremstyle{definition}

  \theoremstyle{remark}

  \theoremstyle{definition} 
  \swapnumbers

\begin{document} 
  \pagestyle{plain}

 \title{Noncommutative orders. A preliminary study} 
 \author{Tomasz Brzezi\'nski}
 \address{ Department of Mathematics, Swansea University, 
  Singleton Park, 
  Swansea SA2 8PP, U.K.} 
  \email{T.Brzezinski@swansea.ac.uk}  
  \date{March 2011}
  \begin{abstract} 
The first steps towards linearisation of partial orders and equivalence relations are described. The definitions of partial orders and equivalence relations (on sets) are formulated in a way that is standard in category theory and that makes the linearisation (almost) automatic. The linearisation is then achieved by replacing sets by coalgebras and the Cartesian product by the tensor product of vector spaces. As a result, definitions of orders and equivalence relations on coalgebras are proposed. These are illustrated by explicit examples that include relations on colagebras spanned by grouplike elements (or linearised sets), the diagonal relation, and an order on a three-dimensional non-cocommutative coalgebra. Although relations on coalgebras are defined for vector spaces, all the definitions are formulated in a way that is immediately applicable to other braided monoidal categories.
  \end{abstract} 
    \maketitle

\section{Introduction}
Mathematical entities of the same kind can be combined to produce a new  object of the same kind: Sets can be combined by the Cartesian product, vector spaces can be combined by the tensor product etc. Objects with a symmetry can also be combined together without the loss of overall symmetry. The Cartesian product of two sets on which a group $G$ acts has a $G$-action provided by the diagonal map.  Two vector spaces on which a Lie group or a Lie algebra is represented can be tensored together to produce a new representation of this group or algebra.\footnote{This is the same principle by which simple quantum mechanical systems can be combined into compound systems without the loss of overall symmetry.}
Formal key property allowing for the latter is the existence of the tensor product or a {\em monoidal structure} in the category of representations of a group (or the category of $G$-sets in the former case). With the birth of quantum groups it has been realised that also their representation spaces can be tensored together. Furthermore, the representations of quantum groups (in the strict sense, that is (dual) quasitriangular Hopf algebras) enjoy a specific symmetry, given in terms of a {\em braiding}. It has been then observed that many noncommutative effects can be explained by using a non-trivial braiding. In other words, the noncommutativity can be incorporated into the notion of symmetry in the category of vector spaces. This has led to  very successful an rich theory of {\em braided groups} and {\em braided geometry} initiated by S.\ Majid \cite{Maj:bra}, \cite[Chapter~10]{Maj:fou}, \cite{Maj:mil}.

Noncommutative geometry can be understood as a {\em linearisation} of the classical (commutative) geometry: the set theoretic notions of classical geometry connected by the Cartesian product are now replaced by (noncommutative) objects belonging to a category with a monoidal or tensor product.\footnote{This is very reminiscent of the transition from classical to quantum mechanics, where classical observables understood as functions on the (symplectic) phase space are replaced by (linear) operators acting on Hilbert spaces in quantum mechanics.}
It is the flexibility of tensor product  that allows for emergence of noncommutative effects. 
Quantum groups or Hopf algebras are a prime example of this  linearisation process. The axioms of a Hopf algebra can be obtained immediately from that of a group by replacing the Cartesian by the tensor product in the axioms of a group written in a fully element-independent way. While there is only one possibility of defining a comultiplication on a group (and this comultiplication or the diagonal map is necessarily cocommutative), the use of the tensor product opens up many new possibilities. 

A programme of developing noncommutative geometry through mo\-noid\-al categories has been initiated recently by T.\ Maszczyk and described in his paper \cite{Mas:non} and numerous lectures. This programme is synthetic in nature, by which we mean that geometric forms are realised in a suitable category \cite{Koc:syn}. With this in mind in \cite{Brz:syn} we have reformulated the notion of a principal bundle within a monoidal category and showen that -- in an appropriately chosen category -- this notion coincides with that of a {\em faithfully flat Hopf-Galois extension} or a {\em principal comodule algebra} which is by now widely accepted as a suitable noncommutative version of a principal fibre bundle (a mathematical object underlying pure gauge theory); see e.g.\  \cite{BrzMaj:gau}, \cite{HajKra:pic}. 

The aim of these notes is by far more elementary. We would like to argue what should be meant by a partial order and by an equivalence relation within a braided monoidal category. A similar problem of formulating generalised or quantum relations in terms of von Neumann (operator) algebras was addressed recently in \cite{Wea:qua}.

The paper is organised as follows. In Section~2 we formulate the standard definitions of equivalence relations and orders on sets in a way that is practiced in category theory; see \cite[Section~2.5]{Bor:han2}. Essentially this is a formulation which avoids using typical set-theoretic concepts such as an element, and replaces them by conditions on functions (or morphisms) and by universal constructions. In Section~3 we translate the abstract formulation of relations on sets to the case of vector spaces. The guiding principle here is that coalgebras should be understood as noncommutative sets. Again, this translation avoids the use of concepts typical for vector spaces, and hence is applicable to any braided monoidal category (which admits some universal constructions). The translation is illustrated by examples in Section~4. Section~5 contains a glossary of coalgebraic terms in hope of making  the material of Sections~3 and 4 more accessible to those, who are not familiar with this terminology.

\section{Orders and equivalence relations on sets}\label{sec.set}
\setcounter{equation}{0}
A relation on a set $X$ is a subset $R$ of the Cartesian product $X\times X$. $R$ is an equivalence relation if it is reflexive, symmetric and transitive, and it is a partial order if it is reflexive, anti-symmetric and transitive. These simple definitions belong to  (any) foundations of mathematics course. In this section we would like to review them from the categorical (and monoidal) point of view; see \cite[Section~2.5]{Bor:han2}.

First fix a singleton set $\{*\}$.  Then any set $X$ comes with two unique mappings: 
\begin{blist}
\item  the {\em diagonal}
$$
\Delta: X\to X\times X, \qquad  x\mapsto (x,x),
$$
\item and the {\em counit}
$$
\eps : X\to \{*\}, \qquad x\mapsto *.
$$
\end{blist}
These two maps can be used to explain the meaning of terms such as reflexive, symmetric, antisymmetric or transitive relation. We can think (more abstractly) about a subset of $X\times X$ as a pair $R$ and an injective (one-to-one) mapping $r: R\to X\times X$. This is how we will understand a relation. The reflexivity of $R$ means that the image of $R$ under $r$ on $X\times X$ contains all pairs $(x,x)$. The latter form the image of the diagonal map. Therefore, to say that the relation $(R,r)$ is {\em reflexive} is the same as to say that there is a mapping $\delta: X\to R$ such that
\begin{equation}\label{reflex}
\Delta = r\circ \delta.
\end{equation}
If $r$ is understood as the inclusion of the subset $R$ into $X\times X$, then $\delta$ is simply the diagonal map.

To explain the meaning of symmetry, recall that $R$ is symmetric if whenever $(x,y)\in r(R)$, then $(y,x)\in r(R)$. The swapping of elements in an ordered pair is encoded in the flip operation
$$
\sigma : X\times X \to X\times X, \qquad (x,y)\mapsto (y,x).
$$
Hence, $(R,r)$ is a {\em symmetric relation} if there exists a mapping $\tau: R\to R$ such that
\begin{equation}\label{sym}
r\circ \tau  = \sigma \circ r.
\end{equation}

The next task is to explain the transitivity. Recall that for this property we need to consider pairs $(x,y)$ and $(y,z)$, such that the second entry of the first coincides with the first entry of the second. If $(x,y)$ and $(y,z)$ are in $r(R)$, then so must be $(x,z)$. To deal with this situation we need to consider {\em pullbacks}. Formally, for any pair of mappings with a common codomain, $f_1: E_1 \to B$, $f_2: E_2 \to B$, their pullback is a 
diagram ({\em i.e.}\ an object $E_1\times_B E_2$ and two morphisms $p_1$, $p_2$ fitting the following diagram)
$$
\xymatrix{ E_1 \times_B E_2 \ar[r]^-{p_1} \ar[d]_{p_2} & E_1 \ar[d]^{f_1} \\
E_2 \ar[r]_{f_2} & B,}
$$ 
with the following universal property. For any set $F$ and mappings $q_1: F\to E_1$ and $q_2: F\to E_2$ such that 
$
f_2 \circ q_2 = f_1 \circ q_1$, there exists a unique mapping $\gamma : F\to E_1 \times_B E_2$ such that $p_2\circ \gamma = q_2$ and $p_1\circ \gamma = q_1$. This situation is usually summarised by the diagram
$$
\xymatrix{ F \ar@{-->}[rd]^-{\exists ! \gamma} \ar@/^/[drr]^-{q_1}\ar@/_/[ddr]_-{q_2} && \\ & E_1 \times_B E_2 \ar[r]^-{p_1} \ar[d]_{p_2} & E_1 \ar[d]^{f_1} \\
& E_2 \ar[r]_{f_2} & B.}
$$ 
Explicitly, for a pair of mappings $f_i: E_i \to B$, $i=1,2$, 
$$
E_1\times_B E_2 = \{(x, y)\in E_1\times  E_2\; |\; f_1(x)=f_2(y)\}.
$$
The projections $p_1$ and $p_2$ are simply restrictions of the canonical projections from the Cartesian product to individual sets. We need to stress, however, that in categories different from the category of sets such an explicit description of a pullback might not be possible.

In the formulation of the transitivity we need to select those pairs of elements of $r(R)$ which have one common element. In other words, if we write
$$
\pi_1 : X\times X\to X, \quad (x,y)\mapsto x,  \qquad \mbox{and}  \qquad \pi_2 : X\times X\to X, \quad (x,y)\mapsto y,
$$ 
then we need to consider all elements of the pullback $R\times_XR$,
$$
\xymatrix{ R \times_X R \ar[r]^-{p_1} \ar[d]_{p_2} & R \ar[d]^{\pi_2\circ r} \\
R \ar[r]_{\pi_1\circ r} & X.}
$$ 
If we write $r(a) = (a_1, a_2), r(b) = (b_1,b_2) \in X\times X$, then
$$
R\times_XR = \{(a,b)\in R\times R \; |\; a_2 = b_1\}.
$$
When viewed inside $X\times X\times X$ the elements of  $R\times_XR$  are triples $(x,y,z)$ such that $x$ is in relation $R$ with $y$ and $y$ is in relation $R$ with $z$. We need to conclude that the outer elements in the triple are also in relation $R$. To pick such elements we define two mappings $r_{L,R}: R\to X$ as the composites
\begin{equation}\label{rLR}
r_L: \xymatrix{ R \ar[r]^-r & X\times X \ar[r]^-{\id\times \eps} & X}, \qquad r_R: \xymatrix{ R \ar[r]^-r & X\times X \ar[r]^-{\eps\times \id} & X}.
\end{equation}
Explicitly, if we write $r(a) = (a_1, a_2)\in X\times X$, then $r_L(a_1, a_2) = (a_1,*) \equiv a_1$ and $r_R(a_1, a_2) = (*,a_2) \equiv a_2$. When restricted to $R\times_ X R$, the Cartesian product $r_L\times r_R$ must give a pair that belongs to $r(R)$. In summary, we say that the relation $(R,r)$ is {\em transitive} provided there exists a mapping $\pi: R\times_XR\to R$, such that
\begin{equation}\label{trans}
r_L\times_Xr_R = r\circ \pi.
\end{equation}
A relation on $X$ understood as a pair $(R,r)$ is an {\em equivalence relation} if it satisfies conditions \eqref{reflex}--\eqref{trans}; see \cite[Section~2.5]{Bor:han2}. 

Recall that anti-symmetry of relation means that whenever $x$ is in relation $R$ with $y$ and $y$ is in relation $R$ with $x$, it must be the case that $x=y$. Thus, to formulate this condition we need to look at elements of $X\times X$ which are both in the image of $r$ and $\sigma\circ r$. We would like to express this property in a way that avoids using elements. This can be done by formulating conditions in terms of mappings as follows.

For any set $Y$ consider mappings $f,g: Y\to R$ such that
\begin{equation}\label{anti1}
r\circ f = \sigma\circ r\circ g.
\end{equation}
 We say that $(R,r)$ is an {\em anti-symmetric relation} if for any pair of mappings $f,g$ satisfying equation \eqref{anti1},
\begin{equation}\label{anti2}
r_L\circ f = r_L\circ g \qquad \mbox{and} \qquad r_R\circ f = r_R\circ g,
\end{equation} 
where $r_L$, $r_R$ are defined in \eqref{rLR}.  Note that the fact that we consider all pairs of mappings $f,g$ satisfying equation \eqref{anti1} amounts to saying that we consider all pairs $(x,y) \in r(R)$ such that also $(y,x) \in r(R)$. 
A relation $(R,r)$ on $X$ that is reflexive, anti-symmetric and transitive is called a {\em (partial) order}.

This seemingly long-winded way of defining  equivalence relations and orders has a few advantages. First, although we talked about relations on sets, we formulated all the axioms in terms of objects, morphisms and some universal constructions. Therefore, this formulation applies to any category with finite limits. Second, and foremost, it can be translated almost {\em verbatim} to any braided monoidal category.

\section{Orders and equivalence relations in monoidal categories}\label{sec.quant} \setcounter{equation}{0}
The aim of this section is to transfer the definitions of an equivalence relation and a partial order to any braided monoidal category. To avoid clattering  the text with abstract (and perhaps not so familiar) notions from category theory, however, we will present the transferred notions in a particular case, namely that of vector spaces over a field $\k$, but in a manner which is applicable to any braided monoidal category (with equalisers). A reader familiar with category theory can easily re-write presented definitions in this generality (assuming that the category in which they are stated has appropriate universal constructions such as equalisers).

\subsection{Relations on coalgebras.} \label{sec.quant1} The category of sets  has the Cartesian product $\times$. This product is associative  and it has an identity (both properties up to bijections) provided by a (fixed) singleton set $\{*\}$. That is, for all sets $X$, $X\times \{*\}\simeq \{*\}\times X\simeq X$. The elements of a pair $(x,y)\in X\times Y$ can be flipped to $(y,x)\in Y\times X$, and this defines a bijective mapping $\sigma: X\times Y \to Y\times X$. In the category of vector spaces there is a tensor product $\otimes$ which serves as a replacement for the Cartesian product. The tensor product is associative and it has an identity $\k$ (both up to linear isomorphisms), {\em i.e.}\ for all vector spaces $V$, $V\ot \k \simeq \k\ot V\simeq V$ by linear isomorphisms. One can also flip tensors by the linear operation $\sigma: V\ot W\to W\ot V$, $v\ot w \mapsto w\ot v$.  The existence of tensor products and the flip (or, more generally, braiding) is all that is needed to define relations in vector spaces. This is essentially a translation of set theoretic notions based on replacing $\times$ by $\otimes$ and $\{*\}$ by $\k$.

We noticed at the beginning of Section~\ref{sec.set} that any set is a coalgebra through the diagonal map in a unique way. It is therefore natural to consider coalgebras $(C,\Delta_C,\eps_C)$ as replacing sets in the category of vector spaces. Note that on a given vector space $C$ one can define various comultiplications, hence we need to specify both $\Delta_C$ and $\eps_C$ as parts of the initial datum. Again in sets any function $r:R\to X\times X$ can be equivalently described as a function assigning to $R$ the graph of $r$, {\it i.e.}, as
$$
\graph{r}: R\to R\times X\times X, \qquad a\mapsto (a, r(a)).
$$
If $X$ and $X\times X$ are understood as coalgebras (in a unique way), then $(R, \graph{r})$ is a right $X\times X$-comodule ($\graph{r}$ is a coaction). Since the diagonal map is cocommutative it is the same as to say that $R$ is an $X$-bicomodule. It seems therefore natural to define a noncommutative or quantum relation over a coalgebra $(C,\Delta_C,\eps_C)$ as a coalgebra (or a ``quantum set'') $(R,\Delta_R,\eps_R)$ that is a $C$-bicomodule with a left coaction $\Rro: R\to C\ot R$ and a right coaction $\roR : R\to R\ot C$. There is, however, slightly less restrictive definition of a quantum relation to which the formulation of orders and equivalences presented in Section~\ref{sec.set} can be transferred. 

Given a coalgebra and a $C$-bicomodule $R$ one can define a $C$-bicomodule map
$$
r = (\id \ot \eps_R \ot \id)\circ (\Rro \ot \id)\circ \roR: R\to C\ot C.
$$
Here $C\ot C$ is a $C$-bicomodule by $\Delta_C\ot \id$ and $\id\ot \Delta_C$. To define such an $r$ there is no need to require $R$ be a coalgebra. By {\em relation on a coalgebra $C$} we will understand a pair $(R,r)$ consisting of a $C$-bicomodule $R$ and a $C$-bicolinear map $r: R\to C\ot C$ (which might be assumed to be a monomorphism to make closer connection with the set-theoretic case). The readers can easily convince themselves that the existence of such an $r$ is equivalent to the existence of a map $\kappa: R\to \k$: Given $r$, define $\kappa= (\eps_C\ot \eps_C)\circ r$; given $\kappa$, define $r= (\id \ot \kappa \ot \id)\circ (\Rro \ot \id)\circ \roR$.

We are now (almost) in position to formulate axioms for quantum orders and equivalences. First, however, we need to look at pullbacks. As already mentioned in Section~\ref{sec.set}, in the category of sets a pullback of $\alpha: E_1 \to B$ and $\beta: E_2 \to B$ is a subset of $E_1\times E_2$ defined by
$$
E_1\times_B E_2 = \{(x,y)\in E_1\times E_2\; |\; \alpha(x) = \beta(y)\}.
$$
Since all sets are coalgebras, and functions are maps of coalgebras, the set $E_1$ is a right $B$-comodule and $E_2$ is a left $B$-comodule with coactions
$$
\lambda_1 = (\id_{E_1}\times \alpha)\circ \Delta_{E_1} : E_1\to E_1\times B, \qquad x\mapsto (x,\alpha(x)),
$$
and 
$$
\lambda_2 = (\beta\times \id_{E_2})\circ \Delta_{E_2} : E_2\to B\times E_2, \qquad y\mapsto (\beta(y),y).
$$
Thus,
$$
E_1\times_B E_2 
= \{(x,y)\in E_1\times E_2\; |\; (x,\alpha(x),y) = (x,\beta(y),y) \} = E_1\coten B E_2,
$$
where $E_1\coten B E_2$ denotes the equaliser of $\lambda_1\times \id_{E_2}$ and $\id_{E_1}\times \lambda_2$, {\it i.e.}, the {\em cotensor product} of comodules. This indicates that pullbacks in a category of sets (or any category with finite limits) should be translated to cotensor products of comodules in a category of vector spaces or any monoidal category. At this point we need to assume that the monoidal category in question has equalisers (and that they are preserved by the tensor product). This assumption is satisfied  for vector spaces.

Let $(R,r)$ be a relation on a coalgebra $C$. Then:
\begin{description}
\item[R1] $(R,r)$ is said to be {\em reflexive} if there exists a map $\delta: C\to R$ such that
\begin{equation}\label{q.reflex}
\Delta_C = r\circ \delta.
\end{equation}.
\vspace{-10pt}
\item[R2] $(R,r)$ is said to be {\em symmetric} if there exists a linear transformation  $\tau: R\to R$ such that
\begin{equation}\label{q.sym}
r\circ \tau  = \sigma \circ r.
\end{equation}
\item[R3] $(R,r)$ is said to be {\em transitive} if  there exists a linear map $\pi: R\Box_CR\to R$, such that
\begin{equation}\label{q.trans}
r_L\Box_Cr_R = r\circ \pi,
\end{equation}
where $r_L = (\id\ot\eps_C)\circ r$ and $r_R = (\eps_C\ot\id)\circ r$.
\item[R4] $(R,r)$ is said to be {\em anti-symmetric} if, for any pair of linear transformations $f,g: V\to R$ such that 
$
r\circ f = \sigma\circ r\circ g,
$
\begin{equation}\label{q.anti2}
r_L\circ f = r_L\circ g \qquad \mbox{and} \qquad r_R\circ f = r_R\circ g. 
\end{equation} 
\end{description}
As in the case of set-theoretic relations, $(R,r)$ that is reflexive, symmetric and transitive is called an {\em equivalence relation on the coalgebra $C$}, while $(R,r)$ that is reflexive, anti-symmetric and transitive is called an {\em order on the  coalgebra $C$}.

We close this section by a few comments on conditions R1--R4. First observe that applying $\id\ot \eps_C$ and $\eps_C\ot \id$ to equation \eqref{q.reflex} one obtains that $r_L\circ\delta = r_R\circ \delta =\id$. Thus the maps $r_L, r_R$ associated to a reflexive relation are necessarily epimorphisms (onto). The above equations mean that they are {\em retractions} (and $\delta$ is their common {\em section}). Next, applying $\id\ot \eps_C$ and $\eps_C\ot \id$ to equation \eqref{q.sym} one finds that the map $\tau$ connects $r_L$ with $r_R$ by $r_R=r_L\circ \tau$ and $r_L =r_R\circ \tau$. Again, the application of $\id\ot \eps_C$ and $\eps_C\ot \id$ to the equality $
r\circ f = \sigma\circ r\circ g,
$ in R4 yields
$$
r_L\circ f = r_R\circ g \qquad \mbox{and} \qquad r_R\circ f = r_L\circ g.
$$
Thus, in particular, every relation in which $r_R = r_L$ is anti-symmetric.

\subsection{Quotients by relations} 
In abstract category theory the process which leads to quotient sets is encoded in terms of {\em coequalisers}. A {\em coequaliser} of morphisms  $q_1, q_2: E\to D$ is an object $C$ together with a morphism $\chi : D\to C$ that {\em coequalises} $q_1$ and $q_2$, that is $\chi\circ q_1 = \chi\circ q_2$, and has the following {\em universal property}: For any morphism $p: D\to B$ that coequalises $q_1$ and $q_2$, there exists a {\bf unique} morphism $q: C\to B$ such that $p = q\circ \chi$. For vector spaces over $\k$, the coequaliser of linear transformations $f,g: V\to W$ is simply the quotient space $W/\im(f-g)$ (the cokernel of the difference $f-g$). The map $\chi : W\to W/\im(f-g)$ is the canonical surjection. 

Let $(R,r)$ be a relation on a coalgebra $C$. By the {\em quotient} $C/(R,r)$ we mean the coequaliser of maps $r_R,r_L: R\to C$. This definition is applicable to any category (it might happen, however, that a given relation does not produce a quotient). In vector spaces
$$
C/(R,r) = C/\im(r_L-r_R).
$$
A reader familiar with coalgebraic techniques will easily find that $\im(r_L-r_R)$ is a {\em coideal} in $C$. This means that, for all $x\in \im(r_L-r_R)$, $\eps_C(x) =0$ and 
$\Delta_C(x) \in C\ot \im(r_L-r_R) \oplus \im(r_L-r_R)\ot C.$ The first statement is obvious as $\eps_C\circ r_L = \eps_C\circ r_R$, a linear map we denoted by $\kappa$ earlier. To prove the second, one needs to use the fact that $r$ is a bicomodule map. Therefore, $C/(R,r)$ is a coalgebra (a ``quantum set'') and the canonical projection $\chi: C\to C/(R,r)$, $c\mapsto [c]$, is a coalgebra map.

It might seem surprising that we define a quotient for any relation not only for an equivalence relation as practiced for sets. This is with complete concord with the categorical approach to the quotients. For any relation $(R,r)$ on a set $X$ one can define maps $r_L$ and $r_R$, and calculate their coequaliser. As a result one obtains the quotient of $X$ by the minimal equivalence relation containing $r(R)$. 

\section{Examples}\label{sec.pcalg}
\setcounter{equation}{0}
\subsection{Linearised sets}\label{sec.lin.set}
Any set $X$ can be made into a cocommutative coalgebra $C$ over $\k$ by defining  $C$ as a vector space with basis $X$, {\em i.e.} $\C = \k X$ and requesting that all elements of $X$ be grouplike, {\em i.e.}, for all $x\in X$,
$$
\Delta_C(x) = x\ot x, \qquad \eps_C(x) =1.
$$
Given  a subset $S\subseteq X\times X$, we define a $C$-bicomodule $R$ to be a vector space with the basis $x\ot y$, for all $(x,y)\in S$ and with the the left and right coactions
$$
\Rro(x\ot y) = x\ot x \ot y, \qquad \roR(x\ot y) = x\ot y\ot y.
$$
The inclusion $S\subseteq X\times X$ extends linearly to a bicomodule map
$$
r: R\to C, \qquad x\ot y \mapsto x\ot y.
$$
Suppose that $S$ is an equivalence relation. Then $S$ contains all pairs $(x,x)$, and hence $x\ot x \in R$, and we can define a map
$$
\delta = \Delta_C: C\to R, \qquad x\mapsto x\ot x.
$$
Obviously $\Delta_C = r\circ \delta$ and $(R,r)$ is a reflexive relation on $C$.  Since $S$ is a symmetric relation on $X$, $(x,y)\in S$ whenever $(y,x)\in S$. Therefore, $v\ot w \in R$ if and only if $w\ot v \in R$. Thus,  $\tau : R\to R$ can be defined as a restriction of $\sigma$ to $R$, and then $\sigma\circ r = r\circ\tau$, so $(R,r)$ is a symmetric relation on $C$.  On the basis $x\ot y$, $(x,y)\in S$, of the vector space $R$, the maps $r_L$ and $r_R$ are
$$
r_L: x\ot y \mapsto x, \qquad r_R: x\ot y \mapsto y.
$$
The basis of the cotensor product $R\Box_C R$ is
$$
x\ot y \ot y \ot z, \qquad (x,y)\in S, (y,z)\in S.
$$
Since $S$ is transitive, also $(x,z)\in S$, hence there is a bicomodule map
$$
\pi: R\Box_C R \to R, \qquad x\ot y \ot y \ot z\mapsto x\ot z,
$$
and
$$
r_L\Box_Cr_R(x\ot y \ot y \ot z) = x\ot z = r\circ \pi(x\ot y \ot y \ot z).
$$
This means that $(R,r)$ is a transitive relation on $C$. Hence $(R,r)$ is an equivalence relation on the coalgebra $C$. The quotient $C/(R,r)$ is the coalgebra spanned by the quotient set $X/S$, {\em i.e.} $C/(R,r) = \k X/S$.

If $S\subseteq X\times X$ is an order on $X$, then corresponding relation $R = \k S$ is a reflexive and transitive relation on $C =\k X$ by the same arguments as before. Since $R$ is a subspace of $C\ot C$, the anti-symmetry of $R$ is equivalent to the requirement that, for every element $w\in R$ such that $\sigma(w)\in R$,
\begin{equation}\label{q.anti.subspace}
r_L(w) = r_L(\sigma(w)), \qquad \mbox{and}\qquad r_R(w) = r_R(\sigma(w)). 
\end{equation}
By the anti-symmetry of $S$, an element $w\in R$ has the property that $\sigma(w) \in R$ if and only if it is flip-invariant, {\em i.e.}\ $w = \sigma(w)$. Then, obviously, equations \eqref{q.anti.subspace} are automatically satisfied. Therefore the linearisation of an order $S$ on $X$ is an order on the coalgebra $\k X$.

\subsection{The diagonal relation}
Let $C$ be a cocommutative coalgebra, and consider the {\em diagonal relation} $(C,\Delta_C)$ on $C$. The diagonal relation is reflexive (take $\delta=\id$) and symmetric (take $\tau = \id$). Both $r_L$ and $r_R$ are identity maps on $C$. Since $C\Box_CC$ is isomorphic with $C$, with the isomorphism given by $\Delta_C$, all elements of $C\Box_CC$ are of the form $\Delta_C(c)$. Choosing $\pi$ as the restriction of $\id\ot\eps_C$ to $\im (\Delta_C) \simeq C\Box_CC$ one immediately sees that the transitivity condition \eqref{q.trans} is satisfied. Hence $(C,\Delta_C)$ is an equivalence relation on $C$. 

In the preceding paragraph we assumed that $C$ is a cocommutative coalgebra. In fact one can easily convince oneself that if $C$ is not cocommutative, then $(C,\Delta_C)$ cannot be an equivalence relation on $C$: As $r_L$ and $r_R$ are identities, the fact that $(C,\Delta_C)$ is symmetric, would imply that $\tau = \id$; see discussion at the end of Section~\ref{sec.quant1}. Hence $\sigma\circ\Delta = \Delta$. On the other hand, even if $C$ is not cocommutative,  $(C,\Delta_C)$ is a reflexive and transitive relation. Since $r_L=r_R$, the discussion at the end of Section~\ref{sec.quant1} affirms that $(C,\Delta_C)$ is an anti-symmetric relation. Therefore $(C,\Delta_C)$ is an order on any coalgebra $C$.

For the diagonal relation $(C,\Delta_C)$, $\im(r_L-r_R) =0$, hence $C/(C,\Delta_C) \simeq C$.

\subsection{An example of a noncommutative order}
Let $C$ be a three-dimensional coalgebra with basis $x,y,z$ and coproduct and counit
$$
\Delta_C(x) = x\ot x, \qquad \Delta_C(y) = x\ot y + y \ot z, \qquad \Delta_C(z) = z\ot z,
$$
$$
\eps_C(x) = \eps_C(z) = 1, \qquad \eps_C(y) = 0.
$$
Consider a five-dimensional subspace $R$ of $C\ot C$ with basis
$$
x\ot x,\ z\ot z,\ x\ot y + y \ot z,\ y\ot x,\ z\ot x.
$$
One easily checks that $R$ is a $C$-bicomodule with left and right coactions given by the restrictions of  $\Delta_C\ot \id$ and $\id\ot\Delta_C$, {\em i.e.}\ $R$ is a sub-bicomodule of $C\ot C$. Take $r$ to be the inclusion map $R\subset C\ot C$. This defines a relation on $C$.

Since $\im \Delta_C\subset R$, the relation $(R,r)$ is reflexive; we need to take $\delta = \Delta_C$. To check the transitivity and anti-symmetry we need to calculate maps $r_L$ and $r_R$. The values of these maps on the specified basis of $R$ are recorded in the following table:\\~

\begin{center}
   \begin{tabular}{c|c|c}
  & $ r_L $ & $r_R$\\
 \hline 
 $x\ot x$ & $x$ & $x$\\ 
 $z\ot z$ & $z$ & $z$ \\
 $x\ot y + y \ot z$ & $y$ & $y$\\
 $y\ot x$ &$y$ & 0\\
 $z\ot x$ & $z$ &$x$
 \end{tabular}
\end{center}~\\~
The space $R\Box_C R$ is six-dimensional with a basis:
$$
x\ot x\ot x\ot x,\ z\ot z\ot z\ot z, \ y\ot x\ot x\ot x,\ z\ot z\ot z\ot x,\ z\ot x\ot x\ot x
$$
and
$$
x\ot x\ot (x\ot y + y \ot z) + (x\ot y + y \ot z)\ot z\ot z.
$$
Applying $r_L\ot r_R$ to the elements of this basis in each case one obtains an element of $R$. Since $r$ is the inclusion, the required map $\pi: R\Box_C R\to R$ is simply the same as the restriction of $r_L\ot r_R$ to $R\Box_C R$. Hence $(R,r)$ is a transitive relation.

All vectors $w$ in $R$ with the property that also $\sigma(w) \in R$ have the form
$$
w = \lambda x\ot x + \mu z\ot z, \qquad \lambda,\mu \in \k.
$$
Since $\sigma(w) =w$,  then obviously $r_L(w) = r_L(\sigma(w))$ and $r_R(w) = r_R(\sigma(w))$, which implies that the relation $(R,r)$ is anti-symmetric; see Section~\ref{sec.lin.set}.

Therefore, $(R,r)$ is an order on the coalgebra $C$.

By inspecting the table of values of maps $r_L$ and $r_R$ one finds that $\im(r_L-r_R)$ is a subspace of $C$ spanned by $y$ and $z-x$. The quotient $C/(R,r)$ is a one-dimensional coalgebra (isomorphic with $\k$) spanned by a group-like element $u$. The canonical projection comes out as:
$$
\chi: C \to C/(R,r), \qquad x\mapsto u, \quad y\mapsto 0, \quad z\mapsto u.
$$
\section{Coalgebraic glossary} 
\setcounter{equation}{0}
A {\em coalgebra} is a vector space $C$ over a field $\k$ equipped with two linear transformations
$$
\Delta_C: C\to C\ot C, \qquad \eps_C : C\to \k,
$$
known as a {\em comultiplication} and a {\em counit}, respectively. These are required to satisfy the following conditions:
$$
(\Delta_C\ot \id)\circ \Delta_C = (\id\ot \Delta_C)\circ \Delta_C, \qquad (\eps_C\ot \id)\circ \Delta_C = (\id\ot \eps_C)\circ \Delta_C = \id.
$$
A coalgebra is said to be {\em (braided) cocommutative}, provided
$$
\sigma\circ \Delta_C = \Delta_C,
$$
were $\sigma$ is a braiding (a flip operator on vector spaces).

A {\em left $C$-comodule} is a vector space $R$ together with a linear transformation $\Rro: R \to C\ot R$ that satisfies the following coassociative and counital laws:
$$
(\Delta_C\ot \id)\circ \Rro = (\id\ot \Rro)\circ \Rro, \qquad (\eps_C\ot \id)\circ \Rro =  \id .
$$

A {\em right $C$-comodule} is a vector space $R$ together with a linear transformation $\roR: R \to R\ot C$ that satisfies the following coassociative and counital laws:
$$
(\roR\ot \id)\circ \roR = (\id\ot \Delta_C)\circ \roR, \qquad  (\id\ot \eps_C)\circ \roR = \id.
$$

A vector space that is both left and right $C$-comodule is called a $C$-bicomodule if the coactions $\Rro$ and $\roR$ satisfy the following compatibility condition:
$$
(\id \ot \roR)\circ \Rro = (\Rro\ot \id )\circ \roR .
$$
The coalgebra $C$ is itself a $C$-bicomodule with both coactions being equal to $\Delta_C$.

Given two left $C$-comodules $(R,\Rro)$ and $(S,\Sro)$ a linear transformation $f: R\to S$ is said to be {\em (left) colinear} or is called a {\em comodule map} if it commutes with the coactions, that is
$$
\Sro\circ f = (\id\ot f)\circ\Rro.
$$
A right colinear map is defined in a similar way. In case of $C$-bicomodules, a {\em bicomodule map} or a {\em bicolinear map} is a linear transformation which is both left- and right colinear.

If $(R,\roR)$ is a right $C$-comodule and $(S,\Sro)$ is a left $C$-comodule then the {\em cotensor product} $R\Box_C S$ of $R$ with $S$ is defined as an equaliser of the maps $\id\ot \Sro\ , \roR\ot \id: R\ot S \to R\ot C\ot S$. That is
$$
R\Box_C S = \{w\in R\ot S \ | \  (\id\ot \Sro)(w) = (\roR\ot \id) (w)\} \subseteq R\ot S.
$$
If $r: R\to R'$ and $s: S\to S'$ are linear maps then $r\Box_C s$ denotes the restriction of $r\ot s$ to $R\Box_C S$.

\section*{Acknowledgments} 
I would like to thank the organisers of the conference on Geometry and Physics in Krak\'ow, September 2010, for creating pleasant and stimulating  atmosphere. Special thanks go to Andrzej Sitarz for very warm hospitality.

 \end{document}